\documentclass[11pt]{article}
	
	\newcommand{\blind}{0}
	
    \makeatletter
    \renewcommand\section{\@startsection {section}{1}{\z@}%
                                       {-3.5ex \@plus -1ex \@minus -.2ex}%
                                       {2.3ex \@plus.2ex}%
                                       {\normalfont\fontfamily{phv}\fontsize{16}{19}\bfseries}}
    \renewcommand\subsection{\@startsection{subsection}{2}{\z@}%
                                         {-3.25ex\@plus -1ex \@minus -.2ex}%
                                         {1.5ex \@plus .2ex}%
                                         {\normalfont\fontfamily{phv}\fontsize{14}{17}\bfseries}}
    \renewcommand\subsubsection{\@startsection{subsubsection}{3}{\z@}%
                                        {-3.25ex\@plus -1ex \@minus -.2ex}%
                                         {1.5ex \@plus .2ex}%
                                         {\normalfont\normalsize\fontfamily{phv}\fontsize{14}{17}\selectfont}}
    \makeatother
	
	\usepackage{amsmath}
	\usepackage{graphicx}
	\usepackage{enumerate}
	\usepackage{natbib} 
	\usepackage{url} 
	
		\usepackage{lineno,hyperref}
		\usepackage[noend]{algpseudocode}
		\usepackage{multirow}
		\usepackage{subfigure}
		\usepackage{wrapfig}
		\usepackage{amsthm}
		\usepackage{amsmath}
		\usepackage{amssymb}
		\usepackage{setspace}
		
		\newtheorem{assumption}{Assumption}
		\newtheorem{lemma}{Lemma}
		\newtheorem{theorem}{Theorem}
		\usepackage{algorithmicx,algorithm}
	
\begin{document}
		
		\def\spacingset#1{\renewcommand{\baselinestretch}%
			{#1}\small\normalsize} \spacingset{1}
		
		\if0\blind
		{
			\title{Enhanced Global Optimization with Parallel Global and Local Structures}
			\author{\normalsize  Haowei Wang $^a$, Songhao Wang$^b$, Qun Meng$^a$, Szu Hui Ng$^a$\\
			$^a$ Department of Industrial Systems Engineering and Management, National  University \\ of Singapore, Singapore 117576\\
             $^b$ Department of Information Systems and Management Engineering, Southern University\\ of Science and Technology, No. 1088 Xueyuan Ave., Shenzhen, China 518000\\ }
			\date{}
			\maketitle
		} \fi
	
		\if1\blind
		{
            \title{Enhanced Global Optimization with Parallel Global and Local Structures}
			\author{Author information is purposely removed for double-blind review}
			\date{}
			\maketitle
		} \fi
	\begin{abstract}
In practice, objective functions of real-time control systems can have multiple local minimums or can dramatically change over the function space, making them hard to optimize.  To efficiently optimize such systems, in this paper, we develop a parallel global optimization framework that combines direct search methods with Bayesian parallel optimization.  It consists of an iterative global and local search that searches broadly through the entire global space for promising regions and then efficiently exploits each local promising region.  We prove the asymptotic convergence properties of the proposed framework and conduct an extensive numerical comparison to illustrate its empirical performance. 
	\end{abstract}
			
	\noindent%
	{\it Keywords:} Parallel global optimization; direct search; Bayesian optimization; global and local search

	\spacingset{1.5} 

\section{INTRODUCTION}
\label{sec:intro}
\sloppy
Significant advances in computing power and the resulting increasing efficiency of simulation models have allowed more simulation models to be applied to real-time control systems: to reflect the complex dynamics of these systems better and to enable more effective control of the systems.  Efficient control of these real-time systems requires effective responses to system changes over time, and the objective function $f(x)$  is usually  continuous and possibly multimodal. A global optimization problem is to find $x^*$ such that $x^*=\arg\min_{x\in\mathbb{X}_{\Omega}} f(x)$, where  $\mathbb{X}_{\Omega} \in \mathbb{R}^d$ is the design space. For stochastic simulation models with random noise, the goal is to find $x^*= \arg\min_{x\in\mathbb{X}_{\Omega}}E(y(x))$, where $y(x)$ is an observation of the objective function corrupted with noise.

When a control system changes dramatically in real-time,  it requires an optimization algorithm that can adapt to the changes quickly. For example, the Safe Sea Traffic Assistant (S2TA)  detects potential conflicts for vessels using an agent-based simulation model \citep{pedrielli2019real}. If a potential risky conflict is detected, the optimizer is called to find an alternative safer trajectory within five minutes. Similar problems also arise in other transportation systems, such as the automated guided vehicle systems where real-time control is required to respond to dynamically changing transportation requests and avoid potential congestion and deadlock in the systems \citep{fazlollahtabar2015methodologies}.

In order to achieve a short response time, the optimization algorithm needs to return solutions quickly. Locally convergent direct optimization algorithms like random search \citep{andradottir1995}, pattern search \citep{torczon,Taddy2009} and COMPASS \citep{hong2006,Hong2007} can be applied given their fast speed. However, when the problems are multimodal, the locally convergent algorithms can get stuck in low-quality solutions \citep{rees2002best}.  A common approach to overcome this is to design restart schemes where the direct search restarts when the search is trapped in sub-optimal areas \citep{marti2003multi,gyorgy2011efficient}. While restarting the local search can better explore and provide a broader search of the function space, the selection of the restart locations is still almost random and largely independent of previously visited locations that have valuable information. Previously visited locations can provide information on non-promising locations (to avoid) and give a global picture of the objective function. Hence, a more intelligent approach to executing the local search can be achieved by learning from previous information.

Advances in computing power has enabled the employment of multiple computing cores to do parallel optimization for faster decision-making.  With the same wall clock time, parallel optimization can search the design space more adequately and facilitate the potential discovery of multiple local solutions to return a better decision. Like the restart scheme mentioned above, in parallel optimization, a critical issue is to design multiple search procedures running on the multiple computing cores, e.g., select multiple starting points for each paralleled local search. An intelligent paralleled search scheme requires efficient use of the previous information in each local core and adequate sharing of information among different cores. 

In this work, we propose a parallel optimization framework that combines fast locally search algorithms and a  Bayesian parallel optimization method. This framework leverages both global and local information from the Additive Global and Local Gaussian Process (AGLGP) model, which stores information of previous observations, to intelligently determine locations to start the paralleled direct search.  The AGLGP model will be refitted with additional observations after an iteration of a global search and a local search, which forms a natural way to communicate between different cores. Such an approach offers a robust and effective alternative to algorithms based only on either a global search or a local search.

There are three main contributions of this work. First, we provide a Parallel Global and Local Optimization (PGLO) framework where a Bayesian parallel Optimization method is used to help intelligently select locations to start the direct search algorithms. This parallel implementation helps to explore the function space more quickly to provide more robust solutions.  
Second, we provide a theoretical convergence guarantee for the proposed framework where any local direct search methods with good empirical performance can be adopted. To test this framework, a specific PGLO with pattern search (PS) is used for an illustration (PGLO-PS).  The third contribution of this paper is an extensive numerical study of PGLO when using one, four, and eight processors on 5 test problems. We first test PGLO-PS empirically and then compare its performance with three metamodel-based optimization algorithms: Two-stage Sequential Optimization (TSSO) \citep{Quan2013}, Expected Quantile Improvement (EQI) \citep{Picheny2013} and Combined Global and Local Search for Optimization (CGLO) \citep{meng2016} with one processor. We further compare it with multi-expected improvement \citep{Ginsbourger} and two multistart implementations of parallel pattern search methods \citep{Hough2001} on 5 test problems with 1,4 and 8 processors. Although we test the examples specifically with pattern search, our parallel framework has theoretical guarantee for a general class of direct search methods.

This work follows the paper \citep{meng2017enhancing} that was published in the Proceedings of the 2017 Winter Simulation Conference. With respect to the conference paper, this manuscript goes beyond in the following main aspects:
\begin{enumerate}
	\item This manuscript formalizes the PGLO algorithm which was only superficially presented in \cite{meng2017enhancing}. Furthermore, \cite{meng2017enhancing} only focused on the enhanced pattern search algorithm. In our current paper we generalize the algorithm so that a wider variety of direct search methods, can be adopted in the proposed PGLO framework. 
	\item This manuscript also provides for the first time the comprehensive convergence proofs and asymptotic properties of the proposed framework. These results apply to the general class of direct search methods (including pattern search).
	\item The empirical analysis is conducted over a wider range of experiments (including higher order examples) with a wider set of competitive algorithms including several serial meta-model bases algorithms, parallel BO algorithms and traditional parallel direct search algorithms.  	
\end{enumerate}

The rest of the paper is organized as follows. Section \ref{sec:basics} discusses the desirable properties of global and local search for global optimization algorithms. The direct search methods and the AGLGP model are also reviewd. Section \ref{sec:pglo} presents the PGLO algorithm and Section \ref{sec:parallel_converge} shows its global convergence. Sections \ref{sec:pglo_num} presents the numerical experiments. Section \ref{sec:conclusions} concludes the whole paper. 

\section{Backgrounds and basics}
\label{sec:basics}
\subsection{Motivating example for combining global and local search}
\label{sec:prop}
Consider the following example for $x\in  [0,1]$:
\begin{equation}
y(x)=(2x+9.96)\cos(13x-0.26)+\epsilon(x). \label{example}
\end{equation}
Here, $\epsilon(x) \sim N(0,4)$ represents the normal noise. The  local and global minimums are at 0.2628 and 0.7460, respectively. Seven points from a Latin Hypercube Design, with 10 replications at each of them, are given as inital design points for illustration.  

\begin{figure}[htbp]
	\centering
	\makebox[\textwidth]{
		\resizebox{1\linewidth}{!}{
			\includegraphics[scale=0.7]{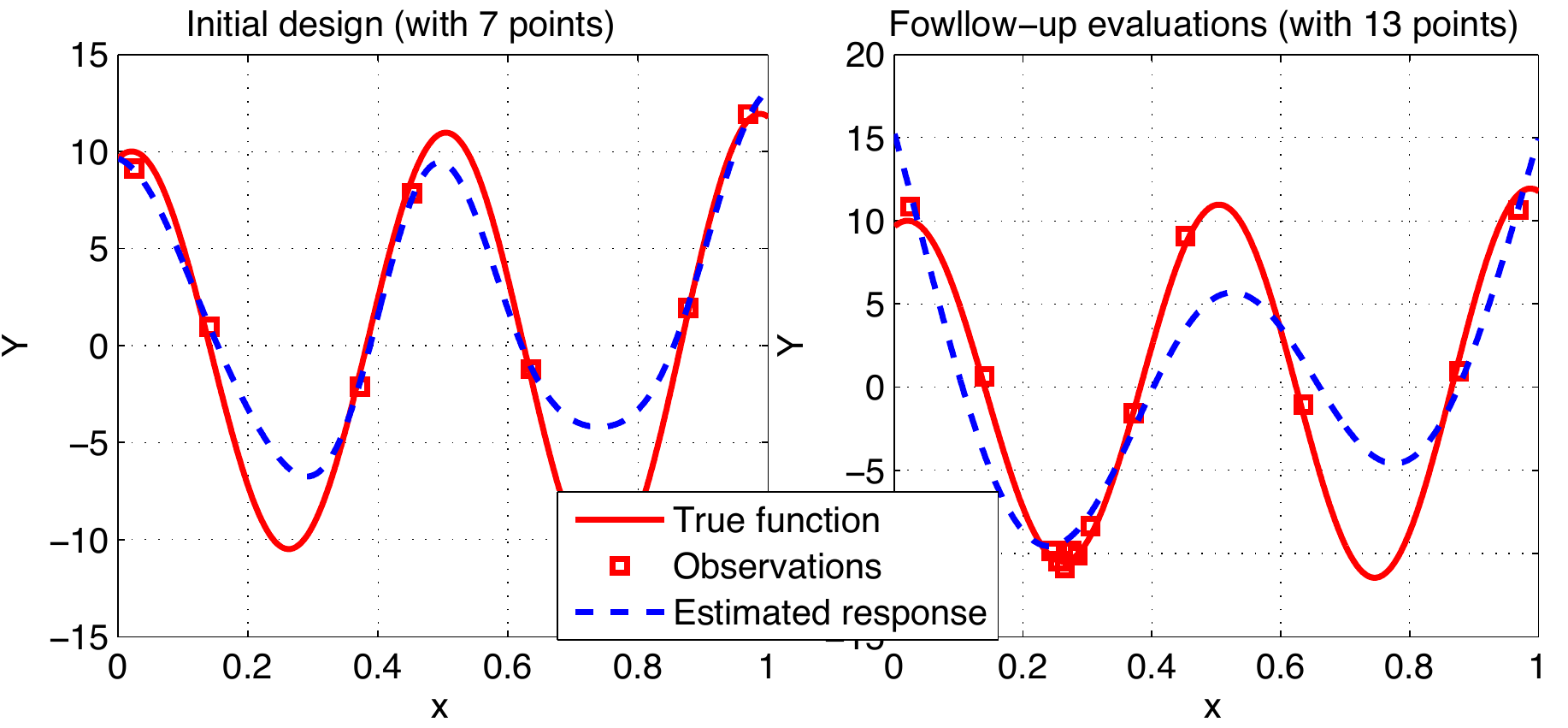}
	}}
	\caption{An illustration of classic pattern search with the seven initial design points and following  6 points chosen by classic pattern search}\label{ps6}
\end{figure}

Fast locally convergent direct search methods are good at quickly identifying and exploiting local optimal regions. Take pattern search \citep{torczon} as an example. It quickly selects points near the current best with a specified lattice for evaluation.
Figure \ref{ps6} shows that the pattern search method selects a sequence of points  exploiting one local optimal region, following the given inital design points.
In this experiment, the point with the smallest sample mean (current best) among the seven initial points is around 0.38, which is closer to the local optimal solution (0.2628) than the true global optimal solution (0.7460).  Starting from this `current best', pattern search gets stuck in the wrong global optimal region.
It can be expected that if multiple pattern search procedures, with additional computing resources, are performed in different local optimal regions simultaneously, the entire space will be explored more sufficiently to find the global optimum. This motivates us to develop a framework to first identify several promising local regions and then use direct search methods to exploit these regions in parallel.

\begin{figure}[htbp]
	\centering
	\makebox[\textwidth]{
		\resizebox{1\linewidth}{!}{
			\includegraphics[scale=0.7]{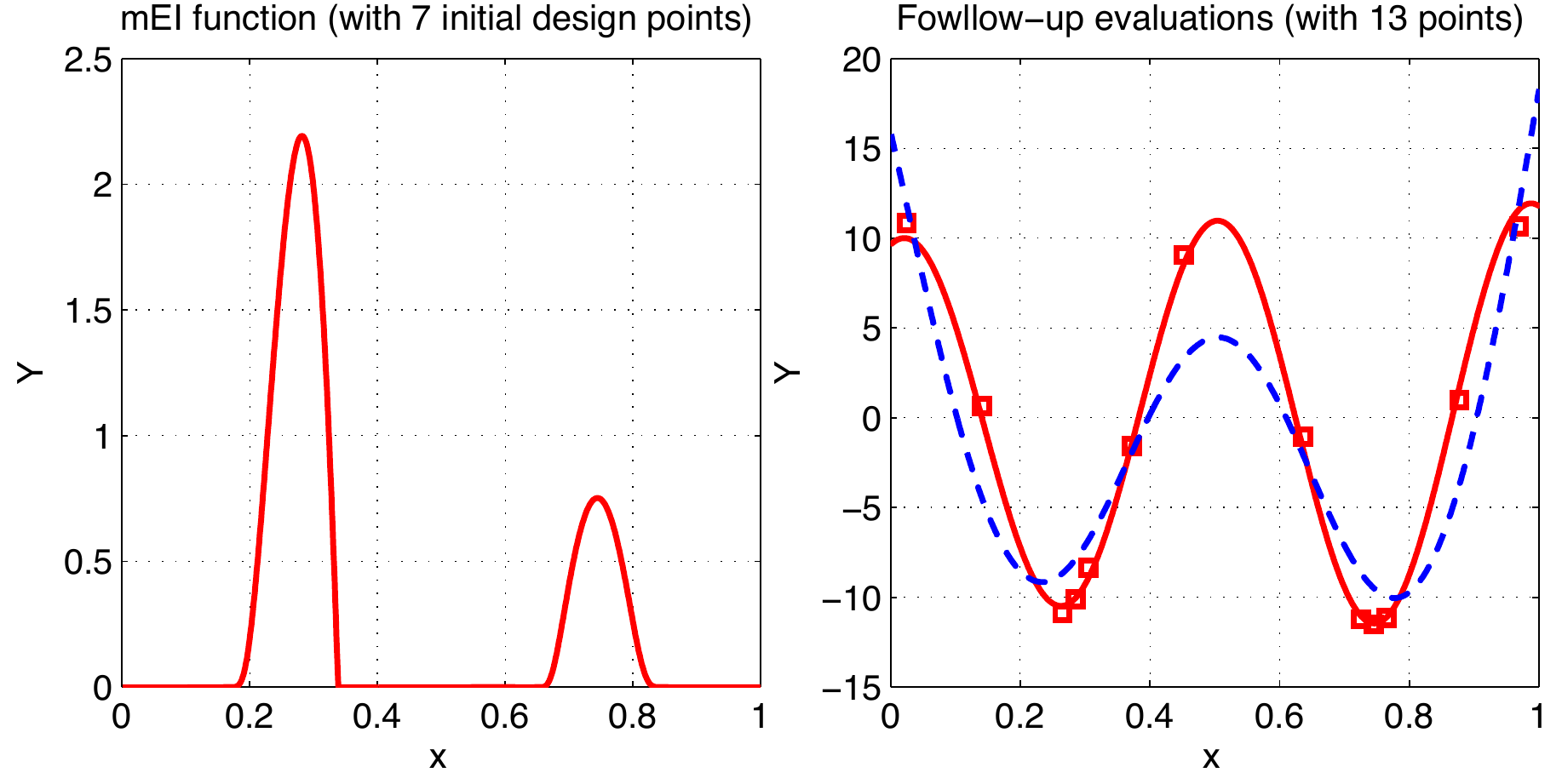}
	}}
	\caption{The mEI criterion with the seven initial design points and following  6 points chosen by parallel pattern search}\label{pps}
\end{figure}

Specifically, we are inspired by a recent Bayesian optimization algorithm, Combined Global and Local Search for Optimization (CGLO) \citep{meng2016}. CGLO employs the AGLGP model to seperately emulate the global and local information of the objective functions (see Section \ref{sec:AGLGP}). CGLO first finds the promising region (i.e., the region with low values or high uncertainty) with the global model and then refines the search with the local model in that region. If several promising regions can be identified with the global model,  the direct search algorithms can be used to search these regions in parallel. 
Take the same example in Equation \eqref{example} for illustration. Left panel of Figure \ref{pps} shows the the modified expected improvement (mEI) function \citep{Quan2013}, the searching criterion used in CGLO, based on the same initial design points as above. Points with larger mEI values are believed to bring more improvement, and thus the point maximizing the mEI function will be selected by CGLO for evaluation. 
Suppose we have two parallel processors and further suppose that each processor exploits the neighborhood of one of the two local optimizers of the mEI criterion. The two series of design points selected by the two parallel pattern search procedures are shown in Figure \ref{pps} right. We observe that the true global optimal region is explored and the optimization process is much improved. Hence, combining pattern search and the global search under the AGLGP model can potentially better explore the whole design space to quickly locate the global optimum. 

Here, we simply choose the two local maximums of the mEI as  starting points for the direct search methods. We will propose a more sophisticated  $q$-points selection criterion, q-mEI, in Section \ref{sec:pglo_initial}. Moreover, the incorporation of pattern search is a specific example for illustration in the development of PGLO. In fact, any general direct search method with good empirical performance can be adopted within our proposed framework. 

\subsection{Direct Search Methods }
\label{sec:ds}
We now present some desirable properties of the direct search algorithms that can be incorporated in PGLO. As discussed in Section \ref{sec:prop}, the direct search methods are expected to exploit local optimal areas more efficiently. Locally convergent random search algorithms, such as random search \citep{andradottir1995}, COMPASS \citep{hong2006}, the adaptive hyperbox algorithm \citep{xu2013}, and the R-SPLINE algorithm \citep{wang2013}, often focus on a local neighborhood to find better solutions. To encourage efficient exploitation, such algorithms need to satisfy the following desirable conditions. Define the \textit{$a$-neighboorhood} of $x\in\mathbb{X}_{\Omega}$ as
\begin{equation}
a\textrm{-}\mathcal{LN}(x)=\{x':x'\in \mathbb{X}_{\Omega}\text{ and }\|x-x'\|\leq a\}. \label{r}
\end{equation}
We then call $x$ a \textit{local minimizer} if $x\in \mathbb{X}_{\Omega}$ and $f(x) \leq f(x')$ for all $x' \in a\textrm{-}\mathcal{LN}(x)$ for some $a>0$. Denote $\mathbb{M}$ as the set of local minimizers in $\mathbb{X}_{\Omega}$. Locally convergent algorithms are guaranteed to converge to solutions in $\mathbb{M}$. \cite{Hong2007} gave a general framework for local convergence and argued that selection of a feasible solution need to satisfy following conditions: for iteration $j$ of the direct search, 

[C1] there exists a deterministic sequence of nested sets, $\mathcal{B}_1\subset\mathcal{B}_2\subset\cdots$, such that the candidate set $\mathcal{C}_j$, where the feasible solution is chosen from, satisfy $\mathcal{C}_j\subset \mathcal{B}_j$, and $\mathbb{X}_{\Omega}\subset\cup_{j=1}^{\infty}\mathcal{B}_j$;

[C2] the sampling distribution $F_j$ guarantees that $Pr\{x\in S_j\}\geq \epsilon $ for all $x \in[ a\textrm{-}\mathcal{LN}(x^*_{j-1})\backslash S_{j-1}]$ for some $\epsilon>0$ that is independent of $j$, where $S_{j-1}$ is the set of solutions estimated through iteration $j-1$.

C1 ensures that the candidate set $\mathcal{C}_j$ is bounded. It prevents sampling points that are arbitrarily far away. C2 ensures that all local neighbors of $x^*_{j-1}$ that have not been evaluated through iteration $j-1$ have a probability of being chosen that is bounded away from zero. See \cite{Hong2007} for detailed discussions.

Algorithms satisfying C1 and C2 have desirable local convergent properties. However, other algorithms without this guaranteed convergence can also have efficient practical performance, as noted in \cite{Hong2007}. These algorithms are still attractive to employ 
in practice and improve practical overall convergence. We note that such algorithms can also be applied in PGLO. With the explorative characteristic of the global model, we show later in Section 4 that PGLO built on direct search methods with unknown local convergence will also converge to the global optimal point.

\subsection{Additive Global and Local Gaussian Process (AGLGP) Model }
\label{sec:AGLGP}
The AGLGP model proposed by \cite{meng2015} can capture both the global and local trends of the responses, providing an effective surrogate model for non-stationary stochastic systems. It models the response as follows:
\begin{gather}
y(x)=f(x)+\epsilon(x)=f_{global}(x)+\sum_{k=1}^Kw_kf^k_{local}(x)+\epsilon(x),\quad w_k=\left\{\begin{array}{rcl} 
1,\quad x\in \mathbf{D}_k \\ 
0,\quad x\notin \mathbf{D}_k
\end{array}\right.\notag\label{model}
\end{gather}
where $f(x)$ is the mean response and $\epsilon(x) \sim N(0,\sigma^2_\epsilon(x))$ is the random noise. Here, $f_{global}(x)$ models the global trend and $f^k_{local}(x), k = 1,\cdots, K$, models the residual process in region $\mathbf{D}_k$, where $\cup_{k=1}^K\mathbf{D}_k=\mathbb{X}_{\Omega}$. AGLGP assumes that $f_{global}(x)$ is a deterministic GP model with mean $\mu$ and covariance $R_g(x_i,x_j) = \sigma^2corr_g(x_i,x_j;\pmb{\theta})$, and that the local model $f^k_{local}(x)$ is a stochastic GP model with 0 mean   and covariance $R_l^k(x_i,x_j) = {\tau_k}^{2}corr_l^k(x_i,x_j;\pmb{\alpha}_k)$. Here,  $\sigma^{2}$ and ${\tau_k}^{2}$ are the variances of $f_{global}(x)$ and  $f^k_{local}(x)$, respectively;  $\pmb{\theta}$ and $\pmb{\alpha}_k$ are the sensitivity parameters of $f_{global}(x)$ and  $f^k_{local}(x)$, respectively. The three different components, $f_{global}(x)$, $f^k_{local}(x)$, and $\epsilon(x)$, are mutually independent.

\begin{figure}[h]
	\centering
	\label{xiong1} 
	\includegraphics[scale=0.6]{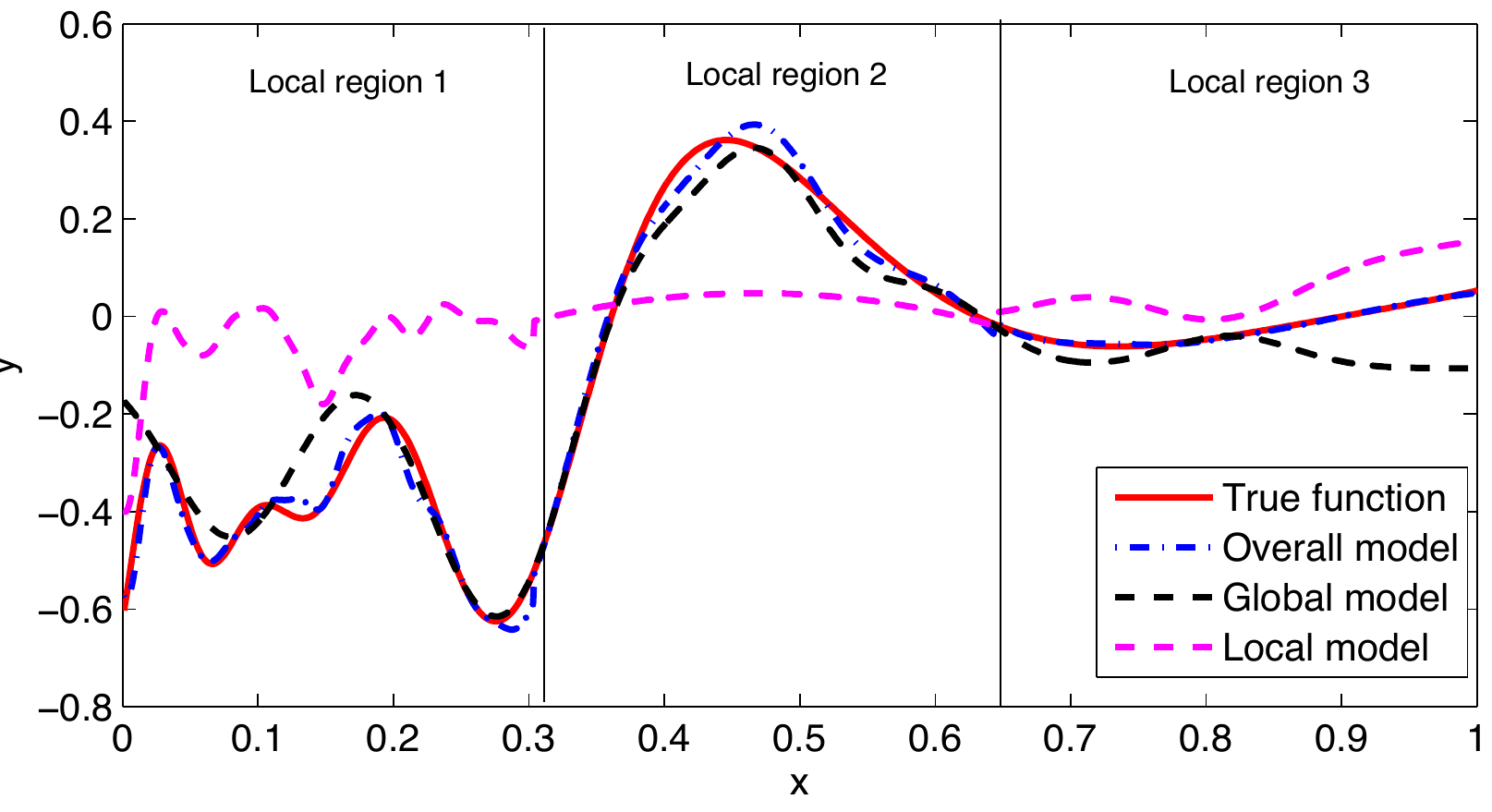}
	\caption{The AGLGP model for  $y(x)=\sin(30(x-0.9)^4)\cos(2(x-0.9))+(x-0.9)/2$ }
	\label{aglgp} 
\end{figure}

Suppose we have $n$ initial design points through Latin Hypercube Sampling over the entire space. $K$ nonoverlapping local regions are selected to segment the entire space via K-means clustering and Support Vector Machine.
The key idea of the AGLGP model is to choose $m$ inducing points first based on the initial design points to estimate the global model, and then model the residuals in each local region with a local model. Readers can refer to the original paper of \cite{meng2015} for details on how to choose the local regions and inducing points. The overall AGLGP predictor $\hat{y}(x)$ for $\forall \,x\in \mathbf{D}_k$ can be derived as follows:
\begin{gather}
\hat{y}(x)=\hat{y}_g(x)+\hat{y}_l(x) \label{overall}\\
\hat{y}_g(x)=\mu+\mathbf{g'Q_m^{-1}G_{mn}(\boldsymbol{\Lambda+\Sigma_{\epsilon}})^{-1}}(\mathbf{y-1'}\mu),\label{gm}\\
\hat{s}_g^2(x)=\sigma^{2}-\mathbf{g'G_m^{-1}g+g'Q_m^{-1}g} \label{gmse}\\
\hat{y}_l(x)=\mathbf{l'(L+\boldsymbol{\Sigma_{\epsilon}})^{-1}(\boldsymbol{\Lambda+\Sigma_{\epsilon}}-\mathbf{G_{nm}Q_m^{-1}G_{mn})}(\boldsymbol{\Lambda+\Sigma_{\epsilon}})^{-1}}(\mathbf{y-1'}\mu), \label{lm}\\
\hat{s}_l^2(x)=\tau_k^{2}-\mathbf{l'(L+\boldsymbol{\Sigma_{\epsilon}})^{-1}l}\label{gmsee}
\end{gather}
where $\mathbf{g}=[R_g(x,x_g^i),\cdots,R_g(x,x_g^m)]'$, $\mathbf{G_m} \in R^{m\times m}$ is a covariance matrix with $(i,j)th$ entry $R_g(x_g^i,x_g^j)$,  $\mathbf{l_k}=[R_l^k(x,x^i_l),...R_l^k(x,x^{N(\mathbf{D}_k)}_l)]'$, $N(\mathbf{D}_k)$ is the number of design points in region $\mathbf{D}_k$,  $n$ is the total number of design points: $n = \sum_{k=1}^KN(\mathbf{D}_k)$. $\mathbf{L_k}$ is the covariance matrix within a local region $\mathbf{D_k}$, $\mathbf{Q_m=G_m+G_{mn}}(\boldsymbol{\Lambda+\Sigma_{\epsilon}})^{-1}\mathbf{G_{nm}}$.  $\mathbf{G_{mn}}\in R^{m\times n}$ is a covariance matrix with $(i,j)th$ entry $R_g(x_g^i,x^j)$, $\mathbf{G_n} \in R^{n\times n} $ is a covariance matrix with $(i,j)th$ entry $R_g(x^i,x^j)$, $\boldsymbol{\Lambda}=diag\{\mathbf{G_n-G_{nm}G_m^{-1}G_{mn}}\}$, $\mathbf{l'}=[\mathbf{0}_{N(\mathbf{D}_1)},...,\mathbf{0}_{N(\mathbf{D}_{k-1})},\mathbf{1}'_{N(\mathbf{D}_{k})},...,\mathbf{0}_{N(\mathbf{D}_{K})}]$ ($\mathbf{0}_{N(\mathbf{D}_k)}$ is a $1\times N(\mathbf{D}_k)$ vector of all 0), $\mathbf{L}=diag\{\mathbf{L_1,\cdots, L_K}\}$ and $\boldsymbol{\Sigma_{\epsilon}}=diag({\sigma}^2_{\epsilon}(x^1),...,{\sigma}^2_{\epsilon}(x^{n}))$. Figure \ref{aglgp} shows an fitting of the AGLGP model for the function $y(x)=\sin(30(x-0.9)^4)\cos(2(x-0.9))+(x-0.9)/2$.


\section{Parallel Global and Local Optimization (PGLO)}
\label{sec:pglo}

In this section, we present the PGLO framework. The initialization of PGLO is to build an initial AGLGP model with a set of $n_0$ initial design points (generated from Latin Hypercube Design). The remaining process of PGLO consists of three stages: a global search stage, a parallel local search stage, and an allocation stage. The global search stage employs the global model \eqref{gm} to find several promising local regions (have either low predicted values or high prediction uncertainty) and distribute the multiple computing cores to them. The local search stage then uses the overall model \eqref{overall} to start each thread of direct search procedure on one computing core. The allocation stage can reduce the noise level to improve the estimation in the previous two steps. 

Suppose that one master and $q$ worker processors are available. Denote $T$ as the total budget (in terms of number of function evaluations) and $n_{max}$ as the budget per iteration.  Table  \ref{PGLOpra} lists the key notations for PGLO and Algorithm  \ref{alg} provides a summary of the proposed framework. In the next subsections, we describe in detail the global search stage (Step 3), the local search stage (Step 4), and the allocation stage (Step 5).

\begin{table}[h]\small
	\centering
	\caption{Key notations for PGLO}\label{PGLOpra}
	\makebox[\textwidth]{
		\resizebox{1\linewidth}{!}{
			\begin{tabular}{c|c}
				\hline
				Inputs & Definition\\
				\hline
				$K$ & the number of local regions\\ 
				$T$ & the total simulation budget \\ 
				$n_0$ & the number of initial evaluation points (Step 1)\\
				$r$ & the number of replications at each initial evaluation point (Step 1)\\
				$n_{max}$ & the budget in each iteration (Step 4)\\ \hline 
				Parameters & Definition \\ \hline 
				$t$ & the index of simulation iterate \\
				$q$	& the number of the worker processors (Step 3)\\
				$k$	& the index of local region (Step 3) \\
				$n_t$ & the number of evaluation points in iteration $t$\\
				$N_t$ & the total number of evaluated points until iteration $t$\\
				
				\hline
			\end{tabular}
	}}
\end{table}

\begin{algorithm}[!t]\small
	\caption{PGLO framework}
	\label{alg}
	\begin{algorithmic}
		\State	\textit{Step 1:} \textit{(Initialization)} Run a Latin hypercube design to obtain $n_0$ initial evaluation points, with $r$ replications at each point. Set iteration $t = 0$ 
		\State	\textit{Step 2:} \textit{(Validation)} Fit an AGLGP model, and use cross-validation (e.g. Leave-one-out cross validation) to ensure a valid AGLGP model.
		\While{Simulation budget is available: set $t = t+1$}
		\State	\parbox[t]{\dimexpr\linewidth-\algorithmicindent}{
			\textit{Step 3:}\textit{(Global Search Stage)} Select $q$ points $x^g_1,\cdots,x^g_q$ based on the multi-point global expected improvement ($q$-$gEI$).  Suppose $q_k$ points among $x^g_1,\cdots,x^g_q$ are located in local region $\mathbf{D_k}$, then $\mathbf{D_k}$ is assigned with $q_k$ processors for the local search. Set $n_t$ = 0.}
		\State	\textit{Step 4:}\textit{(Local Search Stage)} 
		\State \quad \parbox[t]{\dimexpr\linewidth-\algorithmicindent}{
			\textit{4a). (Update Local Models)} Update the local models at $\mathbf{D_k}$ with $q_k>0$ .}
		\State \quad \parbox[t]{\dimexpr\linewidth-\algorithmicindent}{
			\textit{4b). (Select the Starting Points for Direct Search)} In each  $\mathbf{D_k}$, select $q_k$ starting points 
			based on the multi-point modified expected improvement ($q$-$mEI$) criterion. Evaluate these starting points in parallel with $r$ replications at each point. Set $n_t=n_t+q$.}
		\State \quad \parbox[t]{\dimexpr\linewidth-\algorithmicindent}{ \textit{4c). (Select the Follow-up Points by Direct Search)}\\
			\textbf{while} $n_t<n_{max}$, continue the direct search thread on each processor for one further step, i.e., select and evaluate a new point in each of the $q$ processors. Set $n_t=n_t+q$.}
		\State	 \parbox[t]{\dimexpr\linewidth-\algorithmicindent}{\textit{Step 5:}\textit{(Allocation Stage)} Allocate available budget to all the selected points. }
		\EndWhile
		\State	\textit{Step 6:}\textit{(Return)}  Return the evaluated point with best sample mean.
	\end{algorithmic}
\end{algorithm}

\subsection{Global Search Stage}
PGLO evaluates the promisingness of each local region in this stage, based on which the $q$ processes are allocated. Specifically, the optimizer of the multi-point global expected improvement is first determined, i.e., $q$-points are selected by jointly optimizing the searching criterion, and the promisingness of each local region is then quantified by the number of points selected in that region.
Suppose the local region $\mathbf{D_k}$ region has $q_k$ selected points, then $q_k$ processors are allocated to that region. As the global model has smoothed out the local variations, the global search can avoid being stuck in one local region with multiple local optimums.

\subsubsection{Multi-points Global Expected Improvement}
\label{sec:qgei}

To select the $q$ points to explore simultaneously, we propose the following searching criterion called multi-point global expected improvement ($q\textrm{-}gEI$) as an extension of the single-point global expected improvement \citep{meng2016}:
\begin{equation}
\resizebox{0.93\hsize}{!}{$
	q\textrm{-}gEI(x_1^g,...,x_q^g)= \mathbb{E}\{\max\{(y_{gmin}-y_g(x_1^g))^+\cdot\frac{1}{1+e^{n_a(x_1^g)/v-5}},\cdots,(y_{gmin}-y_g(x_q^g))^+\cdot\frac{1}{1+e^{n_a(x_q^g)/v-5}}\}\}.\label{qgei}
	$}
\end{equation}
Here, $y_{gmin}$ is the best global model value at the inducing points, and $y_g(x)$ is normally distributed with mean \eqref{gm} and variance \eqref{gmse}. $(y_{gmin}-y_g(x_i^g))^+$  evaluates the global expected improvement of point $x$. The penalty factor  $\frac{1}{1+e^{n_a(x_i^g)/v-5}}$ gives diminishing returns for points with many existing evaluations around it, where $n_a(x)$ is the number of neighbour design points of $x$ (see Equation \eqref{r}) and $v$ is the penalty parameter. 

It is very computationally expensive  to  maximize the expectation of Equation \eqref{qgei} with $q$ points (of dimension $q\times d$). As a remedy, we select the $q$ points sequentially by optimizing the gEI functions as follows:
\begin{equation*}
\resizebox{1\hsize}{!}{$
	\begin{aligned}
	x^{g}_1=&\arg \max_{x\in \mathbb{X}_{\Omega}} gEI(x)=\arg \max_{x\in \mathbb{X}_{\Omega}}E\{(y_{gmin}-y_g(x))^+\cdot\frac{1}{1+e^{n_a(x)/v-5}}\} ,\\
	x^{g}_{i+1}=&\arg \max_{x\in \mathbb{X}_{\Omega}} gEI(x)=\arg  \max_{x\in \mathbb{X}_{\Omega}} E\{(y_{gmin}-y_g(x))^+\cdot\frac{1}{1+e^{n_a(x)/v-5}}|x^{g}_1,...,x^{g}_{i},y(x^g_1),...,y(x^g_i)\} 
	\end{aligned}
	$}
\end{equation*}

When optimizing the gEI function to obtain $x^g_{i+1}$, the gEI function depends on points $x^g_1, \cdots, x^g_i$, and the responses at these points, which are obviously unavailable yet. In this work, we approximate $y(x^g_i)$ following the similar idea of the \textit{kriging believer} \citep{ginsbourger2010}, i.e. $y(x^g_i) = \hat{y}(x^g_i)$, the AGLGP overall model prediction in Equation \eqref{overall}. It is possible that  the $q$  candidate points  selected are not distinct. To avoid this, if a repeated candidate point is selected, $\tilde{n}_i$ \textit{artificial points} will be added to the penalty factor such that 
\begin{equation}
	gEI(x^g_i)=E\{(y_{gmin}-y_g(x^g_i))^+\}\cdot\frac{1}{1+e^{(n_a(x_i)+\tilde{n}_i))/v-5}}\leq\max_{x\in \mathbb{X}_\Omega \backslash x^g_i} gEI(x).
\end{equation}
Equivalently, this is to assume that there will be additional $\tilde{n}_i$ points selected around $x_i^g$ during the local search stage. This is quite reasonable as at each local search stage, the promising area around $x_i^g$ will be exploited. With a global model that smooths out local variation and a sparse distribution of the global candidates,  the algorithm can search in multiple local regions rather than exploiting in a single region.
\subsection{Parallel Local Search Stage}

This section details the local search stage (Step 4). Recall that in the global search step, $q_k$  processors have been assigned to the local region $\mathbf{D_k}$ and $\sum_k q_k =q$. The larger the $q_k$ is, the more promising the $\mathbf{D_k}$ is. 
In this local search stage, we first select $q_k$ points from  $\mathbf{D_k}$ based on  a multi-point modified expected improvement function ($q\textrm{-}mEI$). Each of these $q_k$ points will be allocated with a computing core and a follow-up direct search is then started from each of the selected points. 

\subsubsection{Starting Points for Direct Search} 
\label{sec:pglo_initial}
The starting points selected should have either lower predicted values or higher prediction uncertainty to facilitate the follow-up direct search to quickly converge to better solutions. Here, we choose the $q_k$ starting points with the following searching criterion, multi-point modified expected improvement (extended from the $mEI$ criterion \citep{Quan2013}):
\begin{align}
q_k\textrm{-}mEI(x_1,...,x_{q_k})= \mathbb{E}\{\max\{(y_{min}^k-z(x_1))^+,...,(y_{min}-z(x_{q_k}))^+\}\}\label{qmei}
\end{align}
where $y_{min}^k$ is the minimum sample mean value in  $\mathbf{D_k}$, and $z(x)$ is  normally distributed  whose mean is Equation \eqref{overall} and variance is $\hat{s}^2_z(x)=\tau_k^2-\mathbf{l'L^{-1}l}$. In order to reduce the computation effort to find a batch $(x_1,...x_{q_k})$,  similar to optimizing the q-gEI in the above, we select the $q_k$ points  by maximizing $mEI(x)$ sequentially as follows:
\begin{equation*}
	\begin{aligned}
	x_1=&\arg \max_{x\in \mathbf{D_k}} mEI(x)=\arg \max_{x\in \mathbf{D_k}}E\{(y_{min}-z(x))^+\},\\
	x_{i+1}=&\arg \max_{x\in \mathbf{D_k}} mEI(x)=\arg  \max_{x\in \mathbf{D_k}} E\{(y_{min}-z(x))^+|x_1,...,x_{i},y_1,...,y_{i}\}.
	\end{aligned}
\end{equation*}

We can observe that the selection of $x_{i+1}$ depends on the selected points $x_1,...,x_i$ as well as their observations $y_1,...,y_{i}$, which are obviously unavailable when we choose the $q_k$ starting points simultaneously. Here, we again adopt the \textit{kriging believer} assumption, i.e. we approximate the `observations' $y_i$ to be $\hat{y}(x_i)$ in Equation \eqref{overall}. This will help avoid selecting new points near existing ones, since the updated variance $\hat{s}^2_z(x)$ (with the artificial observations) at these locations  will be reduced.

\subsubsection{Follow-up Points Selection by Direct search}

After the $q$ starting points are selected by q-mEI, from each of them will start a direct search procedure on one computing core. As a result, each round of these $q$ parallel direct search procedures will select $q$ new points to be evaluated.


Here, we use one of the well-applied direct search methods -- pattern search  (PS) to illustrate our PGLO framework. We detail the PGLO-PS algorithm where PS is used as the direct search method. PS is characterized by its mesh, a lattice where the search is restricted, and polling conditions that specify how the mesh will change. The mesh selection satisfies [C1] and [C2] in Section \ref{sec:ds}. The general steps for PS at iteration $j$ are:

\begin{enumerate}
	\item Generate several candidate points $Q_j$ within a mesh $\mathcal{M}_j$ around the current best point $x_j$.
	\item 
	\begin{enumerate}
		\item Evaluate $Q_j$ with the true models. 
		\item If $\exists x_{j+1}\in {Q}_j$ such that $f(x_{j+1})<f(x_j)$, the search is successful.\\
		Else, the search fails and the mesh is refined, e.g. $\mathcal{M}_{j+1}=\mathcal{M}_j/2$.  Repeat Step (a)
	\end{enumerate}
	\item Update the current best.
\end{enumerate}

The lattices to generate the candidates are predefined and independent of the response  $f$ and existing observed points. The PS can spend most of the computational time on function evaluations as generating trial points is very fast. Therefore, it can exploit a small area efficiently by evaluating a lot of points. However, it has the tendency to get stuck in the local optimal area, especially for multimodel functions. Our PGLO framework can help select the appropriate starting points for PS to improve its efficiency.

Specifically, after selecting $q$ starting points through the q-mEI, we initialize $q$ PS procedures on the selected starting points simultaneously and monitor the progress by recording their mesh sizes, denoted as $\mathcal{M}_{j}^i, i=1,\cdots, q$ at iteration $j$. At the beginning of each pattern search, the starting point is used as the current best point and the evaluation points are generated from  $\mathcal{M}_{j}^i$. The mesh size will shrink if the iteration fails, i.e., no better solution is found. To avoid overly exploiting in unsuccessful areas, the PS stops when
\begin{equation}
\min_{i} \mathcal{M}^i_{j}\leq \mathcal{M}_{min},
\end{equation}
where $\mathcal{M}_{min}$ is a predefined threshold \citep{torczon}. If PS stops when there is remaing budget, i.e., $n_t < n_{max}$, we restart the local search stage. This restart allows the algorithm to jump out of the current local optimums. 

\emph{Remark} It is worth mentioning that all kinds of direct search methods can be adopted under the PGLO framework. Here, we use the pattern search method for illustration.

\subsection{Allocation Stage} \label{allocation}
To reduce the stochastic noise, we use the Optimal Computing Budget Allocation (OCBA) rule \citep{chen2010} to allocate additional simulation replications to all the already evaluated $N_t$ points after the local search stage. Suppose the point $x_i$ has a sample mean $\bar{y}_i$ and a sample variance $\hat{\sigma}^2_{\epsilon}(x_i)$. According to Theorem 1 in \cite{chen2010}, OCBA can optimize the Approximate Probability of Correct Selection (APCS) asymptotically with the following allocation rule:
	\begin{gather}
	\frac{{N}_{OCBC}(x_i)}{{N}_{OCBC}(x_j)}=\left(\frac{\hat{\sigma}_{\epsilon}(x_i)/\Delta_{b,i}}{\hat{\sigma}_{\epsilon}(x_j)/\Delta_{b,j}}\right)^2, i,j\in\{1,2,...,N_t\},\,  i\neq j\neq b\\ 
	{N}_{OCBC}(x_b)=\hat{\sigma}_{\epsilon}(x_b)\sqrt{\sum_{i=1,i\neq b}^n\left(\frac{{N}_{OCBC}(x_i)}{\hat{\sigma}_{\epsilon}(x_i)}\right)^2}\notag\\
	\sum_i {N}_{OCBC}(x_i)=B_{t,a}\notag
	\end{gather}
	where $x_b$ is the best observed point with the best sample meam $\bar{y}_b$, $\Delta_{b,i}=\bar{y}_i-\bar{y}_b$, and ${N}_{OCBC}(x_i)$ is the number of replications at $x_i$ decided by OCBA. $B_{t,a}$ is total budget. OCBA attempts to evaluate more at points with low sample mean or high sample variance as they are more likely to be optimal. 
	
	OCBA was designed to allocate resources to finite alternatives. In PGLO, the design space is assumed to be continuous (with infinite alternatives) and thus an additional assumption should be made to ensure convergence. Denote $M_t(x)$ as the total number of replications received by design point $x$ in the first $t$ iterations. We state the following assumption on the allocation rule:
	\begin{assumption}
		\label{ass:minbdg}
		Suppose the noise variance $\sigma_{\epsilon}^2(x)$ among the input space $\mathbb{X}_{\Omega}$ is bounded, i.e., $\max_{x\in \mathbb{X}_{\Omega}} \sigma_{\epsilon}^2(x)<\infty$. There exists a sequence $\{\kappa_1,...,\kappa_i,...\}$ such that $\kappa_{i+1} \geq \kappa_i$, $\kappa_i \rightarrow \infty$ as $i\rightarrow \infty$ and that $\sum_{i=1}^{\infty}i\exp(-h \cdot \kappa_i) <\infty$, $\forall h>0$.  It follows that $M_t(x)\geq \kappa_{N_t}$ for every of the $N_t$ design points in iteration $t$.
	\end{assumption}
	In the numerical studies (Section \ref{sec:pglo_num}), we select the sequence $\{\kappa_i= 0.05i\}$. To satisfy this assumption, at each evaluation stage, we first ensure that all the $N_t$ evaluated design points receive at least $\kappa_{N_t}$ replications. Then the OCBA is used to allocate the $B_{t,a}$ replications to the $N_t$ design points . These modifications to the budget allocation for each iteration ensures the convergence of the PGLO algorithm.

\section{Convergence of the PGLO Algorithm}
\label{sec:parallel_converge}

In this section, we show that  the PGLO framework is convergent as the number of iterations increases, under some mild conditions. We first list some general assumptions that are used for the analysis in this section.
\begin{assumption} 
	(a) The true response $f$ is bounded.
	
	(b) The variance parameter ($\sigma^2,\tau_k^2$) and lengthscale parameter ($\pmb{\theta},\boldsymbol{\alpha}_k$) for both global and local models are bounded away from zero.
	
	(c) The noise variances $\sigma_{\epsilon}^2(x)$ are known.	
\end{assumption}

Assumptions 2(a) and 2(b) are used in the Efficient Global Optimization (EGO) algorithm to provide convergence guarantee for deterministic response \citep{jones1998efficient}. Assumption 2(c) is to facilitate the analysis. When $\sigma_{\epsilon}^2(x)$ is unknown, the convergence has only been studied empirically to the best of our knowledge \citep{kleijnen2012expected,Pedrielli2018etsso}.

The proof of PGLO framework convergence can be divided into three parts. First, in Lemma 1, we prove that for each local region, if the budget is infinite, the selected design points within the regions are dense.

\begin{lemma}
	Under Assumptions 1 and 2, the sequence of design points within each local region  is dense if this region is given infinitely large budget.
\end{lemma}

Next, we prove that in the global step, each local region will be selected as the promising region infinitely often.

\begin{lemma}
	Under Assumptions 1 and 2, each local region will be selected infinitely often as the promising region to conduct local search. 
\end{lemma}
Combining Lemmas 1 and 2, we see that the sequence of design points will be everywhere dense in the whole design space $\mathbb{X}_{\Omega}$. Finally, in Theorem 1, we prove the convergence of PGLO. 
\begin{theorem}
	Under Assumptions 1 and 2, the optimal value found by PGLO converges to the true global optimum: $\bar{y_t}^*\rightarrow f(x^*)$ w.p.1. as $t\rightarrow \infty$, where $\bar{y_t}^*$ is the current best value and $x^*=\arg\min_{x\in \mathbb{X}_{\Omega}}f(x)$ is the true global optimal solution.
\end{theorem}

Details of the proofs are provided in the appendices.

\section{Numerical Studies} 
\label{sec:pglo_num}

\subsection{Alternative Global Optimization Methods}

In this section, we conduct an extensive comparison with a variety of alternative global optimization methods, including both Gaussian process (GP) based global optimization methods, parallel global optimization methods, and direct search methods.  When evaluating functions, a wait-time of 0.01 seconds is added  to mimic a fast simulation model.

We first compare our proposed PGLO-PS with alternative GP based global optimization algorithms, including TSSO \citep{Quan2013}, EQI \citep{Picheny2013} and CGLO \citep{meng2016}, that are applicable for problems with heteroscadestic case. We then compare its performance with other parallel pattern search methods. Specifically, we include the popular MultPPS-LHS \citep{MatlabOTB}, i.e. multistart parallel pattern search (\textit{MultPPS}) initialized with Latin Hypercube Sample (LHS), and MultPPS-qEI, i.e. \textit{MultPPS} initialized with the q-EI of the AGLGP model. MultPPS-LHS is a direct application of the pattern search method and served as a benchmark algorithm.  MultPPS-qEI is compared to evaluate the benefit of using an iterative of global and local search against a single-level search using q-EI. Besides, Multi-point expected improvement (q-EI) \citep{Ginsbourger} is selected as the alternative parallel global optimization method.

\subsection{Results and Discussions}
To evaluate the performances of different algorithms, we define a reasonable level of accuracy and compute the wall clock time required to find such a solution for each algorithm. Denote $f^*$ as the global optimum and $\hat{f}$ as the estimated optimum, a reasonable solution should has a relative error $|f^*-\hat{f}|/f <1\%$. Also, we use \textit{Speedup} to quantify the enhancement of parallel algorithms. It is defined by $SP=T(1)/T(q)$, where $T(1)$ is the time consumed for a sequential optimization on only one processor and $T(q)$ is the that for that on $q$ parallel processors. We use 30 macro-replications to evaluate the performance. A t-test with a significant level of $\alpha=0.05$ is used to test the significance of difference among different algorithms.

First, we adopt the following example from \cite{sun2014},
\begin{equation}
\max_{0\leq \mathrm{x}_1,\mathrm{x}_2\leq 100}f(\mathrm{x}_1,\mathrm{x}_2)=10\cdot\frac{sin^6(0.05\pi \mathrm{x}_1)}{2^{((\mathrm{x}_1-90)/50)^2}}+10\cdot\frac{sin^6(0.05\pi \mathrm{x}_2)}{2^{((\mathrm{x}_2-90)/50)^2}}.\label{func}
\end{equation}
The global optimum is $f^*(90,90)=20$ and the second best is $f(70,90)=f(90,70)=18.95$. The objective function has 25 local optimal solutions, making it hard for classic optimization methods to find the optimum. The noise term is set to be normally distributed with mean 0 and variance $\sigma^2_{\epsilon}(\mathrm{x}_1,\mathrm{x}_2)=3(1+\mathrm{x}_1/100)^2(1+\mathrm{x}_2/100)^2$.

\begin{table}
	\caption{Average wall clock time (in seconds) to get  a relative error $<1\%$ with one processor}\label{ex:surrogate}
	\centering
	\begin{tabular}{cccc}
		\hline
		PGLO-PS& TSSO & EQI & CGLO\\\hline
		210.36& 1139.13& 1156.52& 252.17\\
		\hline
	\end{tabular}
\end{table}
\begin{table}
	\caption{Average wall clock time (in seconds) to get  a relative error $<1\%$ with different number $q$ of processors}\label{ex:func}
	\centering
	\begin{tabular}{cccccc}
		\hline
		{q}&{PGLO-PS}&{MultPPS-LHS}&{MultPPS-qEI}&q-EI \\\hline
		1&210.36&272.17&252.25 &686.70\\
		
		4&42.43&92.11&53.16&434.92\\
		
		8&37.07&47.93&39.65&213.17\\\hline
		
	\end{tabular}
\end{table}

Table \ref{ex:surrogate} presents the average wall clock time for several GP based global optimization algorithms where the PGLO-PS performs the best with one processor. 
Table \ref{ex:func} presents the average wall clock time for several parallelized algorithms. With only one or four processors, it is difficult for the MultPPS-LHS to explore the entire space as the objective is highly multimodal. Its performance improves a lot with eight processors. The MultPPS-qEI also performs significantly worse than PGLO-PS with one and four processors. With both a global and local search stage, the PGLO-PS can better explore the entire space. The global stage of PGLO-PS helps smooth out the local variations and identify several promising local regions which will be further exploited simultaneously in the local search stage. Therefore, it achieves a balance of exploration and exploitation. In contrast, the MultPPS-qEI with only one search stage could easily get stuck in one of the many local regions. When the number of processors increases to 8, there is no significant difference between the two algorithms as MultPPS-qEI can explore more local regions. The q-EI performs the least efficiently for all cases among different parallelized algorithms. Optimizing the q-EI criterion requires a much longer time, and thus it is more suitable for optimizing functions whose simulation evaluation is very time-consuming. In this work, our focus is on optimizing objective functions that can be quickly evaluated, for which fast direct search methods have the advantage of evaluating more points. In the following example, we omit q-EI as it is not suitable for our problems.


\begin{table}[htbp]
	
	\centering
	\caption{Speedup of parallel optimization algorithm}\label{speedup}
	\begin{tabular}{ccccccc}
		\hline
		q&PGLO-PS&MultPPS-LHS&MultPPS-qEI\\\hline
		4&4.96&2.96&4.75\\
		8&5.68&5.79&6.46\\\hline
	\end{tabular}
\end{table}

Table \ref{speedup} gives a summary of the  $speedup$ with $q=4$ and $q=8$, which measures how well a parallel algorithm scales with additional processors against its sequential counterpart. Both PGLO-PS and MultPPS-qEI have considerable speedup with four processors, but there is a diminishing return with an extra of four processors.
This may because that the two algorithms end up exploits the same optimal area even with additional processors. In contrast, the speedup of MultPPS-LHS with eight processors achieves significant improvement compared to four processors as it has more resources to explore the entire space better.

To comprehensively evaluate the performance of different parallel algorithms, we conduct a more extensive experiment with four commonly used multimodal test functions: Griewank, Ackley, Levy, and Schwefel. All of the four functions have a lot of widespread local minima. Griewank, Ackley, and Levy functions have small local variations and large global variations. The local minima values of the Griewank functions decrease at a regular pace towards its center.  The outer areas of the Ackley are nearly flat, and the function values decrease dramatically around its center. Levy function has a huge variation with respect to one dimension but is very stable with the other. The Schwefel function has both large local and global variations and is more complex compared to the other three functions.  We consider two levels of the noise variance (small and large) for the four test functions based on their global variations.  Specifically, the \textit{small} and noise \textit{large}  variance cases are set to having a maximum of 1\% and 10\%  of the range values of each function respectively. 

\begin{table}[!t]
	\caption{Average wall clock  time (in seconds) and relative speedup to get a relative error $<1\%$  (\textit{small} noise)}\label{ex:optfunc}
	\makebox[\textwidth]{
		\resizebox{1\linewidth}{!}{
			\centering
			
			\begin{tabular}{p{1.4cm}rrrrrrr}
				\hline
				\multirow{2}{*}{Test} &\multirow{2}{*}{$q$}&\multicolumn{2}{c}{PGLO-PS}&\multicolumn{2}{c}{MultPPS-LHS}&\multicolumn{2}{c}{MultPPS-qEI}\\\cline{3-8}
				&& W.C. Time&speedup&W.C. Time&speedup&W.C. Time&speedup\\\hline
				\multirow{3}{*}{Griwank}&1&23.25&-&36.70&-&27.49&-\\
				{}&4&11.20&2.08&30.07&1.22&11.64&2.36\\
				{}&8&10.41&2.23&28.31&1.30&10.54&2.61\\\hline
				\multirow{3}{*}{Ackley}&1&21.06&-&58.98&-&31.92&-\\
				{}&4&12.84&1.64&25.62&2.31&16.27&1.96\\
				{}&8&11.82&1.78&21.58&2.73&12.25&2.61\\\hline
				\multirow{3}{*}{Levy}&1&28.49&-&62.50&-&37.24&-\\
				{}&4&15.54&1.83&37.00&1.69&18.63&2.00\\
				{}&8&15.34&1.86&31.96&1.96&15.71&2.37\\\hline
				\multirow{3}{*}{Schwefel}&1&56.63&-&178.32&-&68.25&-\\
				{}&4&23.65&2.39&57.87&3.08&30.25&2.26\\
				{}&8&18.62&3.04&32.84&5.43&22.14&3.08\\\hline
				
			\end{tabular}
	}}
	%
\end{table}

\begin{table}[!t]
	\caption{Average wall clock  time (in seconds) and relative speedup to get a relative error $<1\%$  (\textit{large} noise)}\label{ex:optfunc2}
	\centering
	\makebox[\textwidth]{
		\resizebox{\linewidth}{!}{
			\begin{tabular}{p{1.05cm}ccccccc}
				\hline
				\multirow{2}{*}{Test} &\multirow{2}{*}{$q$}&\multicolumn{2}{c}{PGLO-PS}&\multicolumn{2}{c}{MultPPS-LHS}&\multicolumn{2}{c}{MultPPS-qEI}\\\cline{3-8}
				&& W.C. Time&speedup&W.C. Time&speedup&W.C. Time&speedup\\\hline
				\multirow{3}{*}{Griwank}&1&33.18&-&56.15&-&32.97&-\\
				{}&4&17.51&1.90&48.71&1.15&23.33&1.41\\
				{}&8&17.18&1.93&49.26&1.14&18.61&1.77\\\hline
				\multirow{3}{*}{Ackley}&1&54.77&-&128.24&-&59.87&-\\
				{}&4&21.43&2.56&66.97&2.31&21.31&2.81\\
				{}&8&16.81&3.26&52.35&2.73&18.40&3.25\\\hline
				\multirow{3}{*}{Levy}&1&72.58&-&150.65&-&53.39&-\\
				{}&4&23.98&3.03&79.51&2.95&33.68&1.58\\
				{}&8&14.56&4.98&53.33&3.87&25.87&2.06\\\hline
				\multirow{3}{*}{Schwefel}&1&110.85&-&298.36&-&148.73&-\\
				{}&4&35.32&3.14&105.25&2.82&35.86&4.15\\
				{}&8&23.57&4.70&48.80&5.38&26.99&5.51\\\hline
				
			\end{tabular}
	}}
	%
\end{table}

Table \ref{ex:optfunc} and \ref{ex:optfunc2}  show the average wall clock time and speedup for the four test functions with small and large noise variances, respectively.  In general, PGLO has the smallest average wall clock time compared to MultPPS-LHS and MultPPS-qEI for all the test functions with different numbers of processors and different levels of noise. This indicates that when the objective functions exhibit high non-stationarity, the PGLO framework that combines the direct search with both global and local metamodels is more efficient than classic paralleled direct search methods and the direct search method guided by only a global metamodel.   It is worth mentioning that the speedup values for PGLO are small for the Griwank, Ackley and Levy functions with four and eight processors for the small noise case. This is because that all the three functions have unique global minimum that is much smaller than other local mimimums. When adopting the AGLGP to smooth out small local variations, it is very easy for the PGLO to identify the global miminum even with only one processor. Focusing on exploiting one promising local region, the MultPPS-qEI and the PGLO-PS have comparable performance in terms of both wall clock time and speedup. For the Schwefel function which has a more complex structure, the PGLO-PS performs significantly better than the other two algorithms. Comparing different noise levels, the wall clock time for large noise case is usually larger than those with small noise, as more evaluations are needed to optimize more noisy functions. The speedup of the large noise case is generally better than that of the small noise case. This may be because that the accuracy of the AGLGP model is improved with additional evaluations, which helps further to better locate the global optimum.

\section{Summary}
\label{sec:conclusions}

We propose a parallel global and local optimization (PGLO) framework in which the global and local information in the AGLGP model is used to guide the direct search method in selecting starting points. The resulting enhanced direct search optimization method can explore the entire space more efficiently and adequately. We theoretically prove the consistency of the proposed algorithm. A comprehensive numerical study is then conducted, and the PGLO performs the best when applied to optimize several classic numerical functions.

%


\bibliographystyle{chicago}
\spacingset{1}
\bibliography{mybibfile}

\begin{appendix}
	
	\section{Proof of Lemma 1}
	
	Let's consider local region $\mathbf{D_{k}}$ and suppose unlimited budget is given to this region. This means that $\mathbf{D_{k}}$ has been selected as one promising region infinitely often and thus for infinite number of iterations, we will perform local search procedures in this region. Recall that to start the local search procedures, PGLO will first select one or several starting points of the local search (depending on the number of processors assigned to this region from the global search step). To select the starting point(s), PGLO uses multi-point mEI criterion, and if several starting points should be selected, the kirging believer strategy will be adopted. Therefore, the first starting point will always be selected by a traditional mEI criterion applied to this region and this point will be evaluated in the following local search procedure starting from itself. Let us then only focus on this single `first' starting point. In each iteration when $\mathbf{D_{k}}$ is selected as promising region, PGLO will evaluate a point selected by the original mEI criterion (and of course several other points from the following procedures in the local search steps). As $\mathbf{D_{k}}$ will be selected as promising region infinitely often, there are infinite points that are selected by the original mEI criterion. As long as the series of `first mEI points' are dense, all the design points selected in the region will be dense.
	
	The following proof will be divided into three parts following that of Lemma 1 in \cite{meng2016} and that of Theorem 1 in \cite{locatelli1997bayesian}. First, Section A.1 gives an upper boundthat of the mEI value at any unsampled point in $\mathbf{D_{k}}$. Section A.2 then constructs a region close to a evaluated point where all the unobserved points have bounded mEI values. Section A.3 finally finished the proof.
	
	\subsection{The upper bound for mEI value at any unsampled point}
	The mEI at $x$ is:
	\begin{equation}
	\label{detailmEI}
	mEI(x)=E(\max[y_{min}-z(x)],0)= \widehat{s}_z(x) \phi (\frac{y_{min}-\widehat{\mu}(x)}{\widehat{s}_z(x)})+ (y_{min}-\widehat{\mu}(x))\Phi(\frac{y_{min}-\widehat{\mu}(x)}{\widehat{s}_z(x)}),\tag{A.0}
	\end{equation}
	where $y_{min}$ is the current best in $\mathbf{D_{k}}$, $\hat{s}^2_z(x)=\tau^2_{k}-\mathbf{l'L^{-1}l}$ is the local spatial uncertainty, $\widehat{\mu}(x) = \hat{y}_g(x)+\hat{y}_l(x)$, if $\underline{M}\leq \hat{y}_g(x)+\hat{y}_l(x)\leq  \overline{M} $; $\widehat{\mu}(x)=\underline{M}$ if $\hat{y}_g(x)+\hat{y}_l(x) <\underline{M}$; $\widehat{\mu}(x)=\overline{M}$ if $\hat{y}_g(x)+\hat{y}_l(x) >\overline{M}$. This is the overall predictor constrained in $[\underline{M},  \overline{M} ]$ (to avoid severe under or over estimation of the predictor, see details in \cite{meng2016}).
	
	According to Assumption 2(a), the response is bounded. We then select a large value $M$ such that both observations $\mathbf{y}$ and prediction $\widehat{\mu}(x)$, are bounded in $(-M,M)$ for all $x$. Therefore, we see 
	\begin{equation}
	-2M<y_{min}-\widehat{\mu}(x)<2M.\tag{A.1}\label{B1}
	\end{equation} 
	
	On the other hand, denote $x_0$ as the nearest design point in $\mathbf{D_{k}}$ to $x$. We can find that
	$$\hat{s}^2_z(x)=\tau^2_{k}-\mathbf{l'L^{-1}l}<\tau^2_{k}-\tau^2_{k}corr_l^k(x,x_0).$$ Define $\hat{s}^2_0:=\tau^2_{k}-\tau^2_{k}corr_l^k(x,x_0)$. Note that $mEI(x)$ is an increasing function of $(y_{min}-\widehat{\mu}(x))$ and $\hat{s}_z(x)$. It follows that,
	$$mEI(x)\leq \widehat{s}_0(x) \phi (\frac{2M}{\widehat{s}_0(x)})+ 2M\Phi(\frac{2M}{\widehat{s}_0(x)}):= U(x;x_0).$$
	As a result, we find an upper bound for $mEI(x)$ considering $x_0$. 
	
	\subsection{Local region centered at design points with bounded mEI values}
	As $\hat{s}^2_0$ is an increasing function with respect to the distance between $x$ and $x_0$, and that $U(x;x_0)$ increases with $\hat{s}^2_0$, we see that $U(x;x_0)$ increases as the distance between $x$ and $x_0$ increases. Hence, we can build a region around $x_0$ such that within this region, $U(x;x_0)$ values for all unsampled points are smaller than a threshold $c$:
	$$R(x_0,c):=\{x\in \mathbf{D_{k}}|U(x;x_0)<c\} . $$ Therefore, we find regions centered at any design points in $ \mathbf{D_{k}}$ such that the mEI values within are smaller than $c$. With $R(x_0,c)$, we can apply Lemma 1 and Theorem 1 from \cite{locatelli1997bayesian} that if the algorithm stops when the maximum of mEI is smaller than predefined threshold $c$, the points within $R(x_0,c)$ for all $x_0$ will never be selected, and thus the algorithm will terminate with a finite iteration. 
	
	\subsection{Proof of density}
	Next, we prove that for any point $x_a$ in $\mathbf{D_{k}}$ and any $\epsilon>0$, there exist a large value $K$ such that at least one design point $x_1$ is selected before iteration $K$ in region $S=\{x||x_a-x_1|<\epsilon\}$. If before $K$, at one design point have already be selected, the condition holds. Otherwise, according to \eqref{B1}, we can find an lower bound for mEI value at $x_a$:
	$$mEI(x_a)\geq \widehat{s}_1(x_a) \phi (\frac{-2M}{\widehat{s}_1(x_a)})- 2M\Phi(\frac{-2M}{\widehat{s}_1(x_a)}):= P(x_a), $$
	where $\widehat{s}_1^2(x_a)>0$ is the predictive variance if all points in $\mathbf{D_{k}} \setminus S$ are selected as design points, which serves as a lower bound for $\widehat{s}_z^2(x_a)$. Moreover, it is easy to check that $P(x_a)>0$ given $\widehat{s}_1^2(x_a)>0$. In this case, by setting $c<P(x_a)$ and selecting large $K$, within finite iterations $K-1$, the mEI values at all points in $\mathbf{D_{k}} \setminus S$ will be smaller than $c$. Therefore, the next design point must be chosen from $S$. It follows that,
	\begin{equation}
	\mathbb{P}[\min_{x_1\in D}|x_1-x_a|>\epsilon\ \  i.o.]=0, \tag{A.2}\label{B2}
	\end{equation}
	where $D$ is the design set and $i.o.$ is shorthand for infinitely often. With Equation \eqref{B2}, we see that design points are dense almost surely.

	\section{Proof of Lemma 2}
	
	In each global search step of PGLO, multiple global candidate points will be selected through q-gEI criterion. This q-gEI criterion is an extension of the gEI criterion proposed in CGLO \citep{meng2016} to the cases of multiple candidate points. This extension essentially employs the kriging believer technique. It will select the first candidate point using the original gEI criterion and then treat the prediction from the global model as the true global observation at this candidate. To select the following candidate points, the first candidate will be treated as a known point. Therefore, we see that in each iteration of PGLO, there is a candidate point selected by the original gEI criterion and thus the region this point belongs to will be selected as a promising region for local search. We also note that there can be other regions chosen as the promising regions in the same iteration with the remaining candidate points. In short, with unlimited budget, in each iteration, at least one local region is selected using the original gEI criterion. We next consider this 'first gEI region' in each iteration and we prove that each local region will be selected as this first gEI region in unlimited number of iteration. The remianing proof of this lemma will follow that of Lemma 2 in \cite{meng2016} and we will repeat it here.
	
	With a limited number of regions and an infinitely large budget, there must exist some local regions that are chosen as promising regions with infinite times. Otherwise, the algorithm will terminate once each local region is chosen finite times. We prove this lemma by contradiction. Assume there are only two local regions $\mathbf{D_{1}}$ and $\mathbf{D_{2}}$ where $\mathbf{D_{1}}$ is chosen infinite times and $\mathbf{D_{2}}$ is chosen as promising region $K_1$ times (for more than 2 regions, $\mathbf{D_{1}}$ is the union of regions selected infinitely and $\mathbf{D_{2}}$ is the union of regions selected finite times). 
	
	First, we prove that after chosen $K_1$ times, the gEI value at any points in $\mathbf{D_{2}}\cap \mathbb{X}_{G}$ has an lower bound. We rewrite the gEI function as follows:
	$$gEI(x)= E\{\max(y_{gmin}-y_g(x)),0 \}\cdot \frac{1}{1+e^{n_i/v-5}} := gEI_v(x) \cdot gEI_p(x),$$
	where $gEI_v(x)$ is obtained from the EI function and $gEI_p(x)$ is the penalty term. Note that as $\mathbf{D_{2}}\cap\mathbb{X}_{G}$ is no longer selected after $K_1$ times, the $gEI_p(x)$ value will not be further updated. Define $g_p:=\min_{x\in \mathbf{D_{2}}} gEI_p(x)$.
	
	In addition, we see that
	$$\hat{s}_g^2(x_0)=\sigma^{2}-\mathbf{g'G_m^{-1}g+g'Q_m^{-1}g}.$$
	The first two terms $\sigma^{2}-\mathbf{g'G_m^{-1}g}$ have a positive lower bound $\hat{s}_{g,1}^2(x)$. This lower bound is attained when all the points in $\mathbf{D_{1}}\cap\mathbb{X}_{G}$ in the following iterations are selected as inducing points. For the third term, we prove it is positive since $\mathbf{Q_m}$ is positive definite. Recall that,
	$$\mathbf{Q_m=G_m+G_{mn}}(\boldsymbol{\Lambda+\Sigma_{\epsilon}})^{-1}\mathbf{G_{nm}} . $$
	In this equation, $\mathbf{G_m}$ is positive definite as it is a covariance matrix. In addition, $(\boldsymbol{\Lambda+\Sigma_{\epsilon}})$ is a diagonal matrix whose diagonal entries are all positive, and thus, $ \mathbf{pG_{mn}}(\boldsymbol{\Lambda+\Sigma_{\epsilon}})^{-1}\mathbf{G_{nm}p'}>0$ for any vector $\mathbf{p}$. It follows that $\mathbf{Q_m}$ is positive definite. Therefore, we find a positive lower bound $\hat{s}_{g,1}^2(x)$ for $\hat{s}_{g}^2(x)$. 
	
	With the lower bound, following the same reasoning with Appendix A, we find the following positive lower bound for $gEI(x)$ for an arbitrary point $x$ in $\mathbf{D_{2}}\cap\mathbb{X}_{G}$:
	$$gEI(x)> \{\widehat{s}_{g,1}(x) \phi (\frac{-2M}{\widehat{s}_{g,1}(x)})- 2M\Phi(\frac{-2M}{\widehat{s}_{g,1}(x)})  \}\cdot g_p :=T(x).$$
	
	Let's next consider the gEI value for any point $x'\in \mathbf{D_{1}}\cap\mathbb{X}_{G}$. As $\mathbf{D_{1}}$ is selected as promising region infinitely, according to Lemma 1, the design points within are dense. Specifically, we can select a large number of $K_2$ such that for iteration $t>K_2$, for any $x'\in \mathbf{D_{1}}\cap\mathbb{X}_{G}$, there will be at least $N$ design points in $\mathbb{B}(x')$ (recall that $\mathbb{B}(x')$ denotes the neighbors set of $x'$, see Section 4.1) such that
	$$gEI(x') = gEI_v(x) \cdot  \frac{1}{1+e^{N/v-5}}<\sigma^2  \frac{1}{1+e^{N/v-5}}<T(x)  ,$$ 
	where $x$ is an arbitrary point in $\mathbf{D_{2}}\cap\mathbb{X}_{G}$. It follows that, $\max_{x \in\mathbf{D_{2}}\cap\mathbb{X}_{G}}gEI(x)>\max_{x' \in \mathbf{D_{1}}\cap\mathbb{X}_{G}}gEI(x')$. Therefore, in the next iteration, $\mathbf{D_{2}}$ will be chosen as promising region, which contradicts the assumption that $\mathbf{D_{2}}$ will never be selected again and finishes the proof.
	
	\section{Proof of Theorem 1}
	
	With Lemma 1 and Lemma 2, we observe that the design points in the entire design space will be dense asymptotically. The proof of the theorem will then follow that of Theorem 1 in CGLO \citep{meng2016}. We will repeat the proof as follows.
	
	In iteration $t$, there are $N_t$ design points in total (we denote the whole design point set as $D_t$). Denote $\hat{x}_t$ as the current best solution such that $\bar{y}^*_t=\bar{y}_t(\hat{x}_t)=\min_{x\in D_t}\bar{y}_t(x)$. In addition, denote $x_0^t$ as the true best in $D_t$ such that $f(x_0^t)=\min_{x\in D_t}f(x)$. This theorem states that $\bar{y}^*_t\rightarrow f(x^*) $ w.p.1 as $t\rightarrow \infty$, where $x^*=\arg\min_{x\in \mathbb{X}_{\Omega}}f(x)$. We fist prove $\bar{y}^*_t-f(x_0^t)\rightarrow  0 $ w.p.1 in Section C.1 and then $f(x_0^t)- f(x^*)\rightarrow  0 $ w.p.1 in Section C.2. Finally in Section C.3, we finish the proof.
	\subsection{Proof that $\bar{y}^*_t-f(x_0^t)\rightarrow 0 $ w.p.1 as $t\rightarrow \infty$}
	The sufficient condition is that $ \sum_{t=1}^{\infty}\mathbb{P}[|\bar{y}^*_t-f(x_0^t)|>\delta]<\infty$, for any $\delta>0$. Note that
	\begin{equation}
	\begin{split}
	\mathbb{P}[|\bar{y}^*_t-f(x_0^t)|>\delta ] & = \mathbb{P}[|\bar{y}_t(\hat{x}_t)- f(\hat{x}_t)+f(\hat{x}_t)- f(x_0^t)|>\delta ] \\& <  \mathbb{P}[|\bar{y}_t(\hat{x}_t)- f(\hat{x}_t)|>\frac{\delta}{2}] + \mathbb{P}[|f(\hat{x}_t)- f(x_0^t)|>\frac{\delta}{2}]\end{split} \tag{C.1}\label{D1}
	\end{equation}
	
	We first bound the first term in \eqref{D1}. Define $\sigma^2_0:= \max_{x \in \mathbb{X}_{\Omega}}  \sigma^2_\epsilon(x)$. For any $x \in D_t$, $\bar{y}_t(x)- f(x)$ is a Gaussian variable $\mathcal{N}(0, \sigma_\epsilon^2(x)/M_t(x))$, and thus (recall $N_t$ is the total number of design points in iteration $t$)
	\begin{equation*}
	\mathbb{P}[|\bar{y}_t(x)- f(x)|>\frac{\delta}{2}]\leq 2\exp(-\frac{\delta^2 M_t(x)}{8\sigma_\epsilon^2(x) } ) \leq 2\exp(-\frac{\delta^2 \kappa_{N_t}}{8\sigma_0^2 } ).
	\end{equation*}
	Therefore, 
	\begin{equation}
	\mathbb{P}[\max_{x \in D_t}|\bar{y}_t(x)- f(x)|>\frac{\delta}{2}] \leq \sum_{i=1}^{N_t}\mathbb{P}[|\bar{y}_t(x_i)- f(x_i)|>\frac{\delta}{2}] \leq 2N_t\exp(-\frac{\delta^2 \kappa_{N_t}}{8\sigma_0^2 } ) .\tag{C.2} \label{D2}
	\end{equation}
	Hence, 
	$$\mathbb{P}[|\bar{y}_t(\hat{x}_t)- f(\hat{x}_t)|>\frac{\delta}{2}] <\mathbb{P}[\max_{x \in D_t}|\bar{y}_t(x)- f(x)|>\frac{\delta}{2}]\leq  2N_t\exp(-\frac{\delta^2 \kappa_{N_t}}{8\sigma_0^2 } ). $$ 
	
	To bound the second term in \eqref{D1}, define sets $A:=\{|\bar{y}_t(\hat{x}_t)- f(\hat{x}_t)|<\frac{\delta}{5} \}$ and $B:=\{|\bar{y}_t(x_0^t)- f(x_0^t)|<\frac{\delta}{5} \}$. Thus,
	\begin{equation}\mathbb{P}[|f(\hat{x}_t)- f(x_0^t)|>\frac{\delta}{2}]= \mathbb{P}[\{|f(\hat{x}_t)- f(x_0^t)|>\frac{\delta}{2}\} \bigcap \{A\bigcap B \}  ] +\mathbb{P}[\{|f(\hat{x}_t)- f(x_0^t)|>\frac{\delta}{2}\} \bigcap \{A\bigcap B \}^\complement  ].  \tag{C.3} \label{D3}\end{equation}
	We prove the first term in \eqref{D3} is zero by contradiction. When $f(\hat{x}_t)- f(x_0^t)>\frac{\delta}{2}$, $|\bar{y}_t(\hat{x}_t)- f(\hat{x}_t)|<\frac{\delta}{5}$ (set A) and $|\bar{y}_t(x_0^t)- f(x_0^t)|<\frac{\delta}{5}$ (set B), it must be that $\bar{y}_t(\hat{x}_t)>\bar{y}_t(x_0^t)$. This contradict that $\bar{y}_t(\hat{x}_t)$ is the current best. Therefore the first term in \eqref{D3} is zero. Consider the second term in \eqref{D3}, 
	$$\mathbb{P}[\{|f(\hat{x}_t)- f(x_0^t)|>\frac{\delta}{2}\} \bigcap \{A\bigcap B \}^\complement ] < \mathbb{P} [ \{A\bigcap B \}^\complement]<2-\mathbb{P}[A]- \mathbb{P}[B]<4N_t\exp(-\frac{\delta^2 \kappa_{N_t}}{50\sigma_0^2 } ). $$
	The last inequality results from the same reasoning with \eqref{D2} that $ 1-\mathbb{P}[A]= \mathbb{P}[|\bar{y}_t(\hat{x}_t)- f(\hat{x}_t)|>\frac{\delta}{5}]<2N_t\exp(-\frac{\delta^2 \kappa_{N_t}}{50\sigma_0^2 } ) $ (similarly for $ 1-\mathbb{P}[B]$). Therefore, we find a bound for the second term in \eqref{D1}:
	$$ \mathbb{P}[|f(\hat{x}_t)- f(x_0^t)|>\frac{\delta}{2}]<4N_t\exp(-\frac{\delta^2 \kappa_{N_t}}{50\sigma_0^2 } ). $$
	Hence, $$\mathbb{P}[|\bar{y}^*_t-f(x_0^t)|>\delta ]<6N_t\exp(-\frac{\delta^2 \kappa_{N_t}}{50\sigma_0^2 } ) $$
	
	$$\mathbb{P}[|\bar{y}^*_t-f(x_0^t)|>\delta ]< 2N_t\exp(-\frac{\delta^2 \kappa_{N_t}}{8\sigma_0^2 } )+ 4N_t\exp(-\frac{\delta^2 \kappa_{N_t}}{50\sigma_0^2 } ) <6N_t\exp(-\frac{\delta^2 \kappa_{N_t}}{50\sigma_0^2 } ) $$
	
	As a result, by Assumption 1, 
	$$\sum_{t=1}^{\infty}\mathbb{P}[|\bar{y}^*_t-f(x_0^t)|>\delta]<6 \sum_{t=1}^{\infty}  N_t\exp(-\frac{\delta^2 \kappa_{N_t}}{50\sigma_0^2 } ) < \infty. $$
	
	\subsection{Proof that $f(x_0^t)- f(x^*)\rightarrow  0 $ w.p.1 as $t\rightarrow \infty$}
	
	Suppose that the global optimal solution $x^*$ lies in region $\mathbf{D_{k}}$. For any $\delta>0$, we can select a region $S\subset\mathbf{D_{k}}$, such that for all $x\in S$, $|f(x)- f(x^*)|\leq \delta$ and that $|x-x^*|<\eta$ (assuming the response $f$ is continuous). 
	
	According to Lemma 1, for any $\epsilon>0$, we can select a large number $K$ such that there exists at least one design point within region $S$ in the first $K$ selected design points in $\mathbf{D_{k}}$. By Lemma 2, we can find a large number $K_2$, such that in the first $K_2$ iterations of CGLO, the region $\mathbf{D_{k}}$ will be selected as promising region at least $K$ times. Therefore, at least one design point in $S$  can be selected before iteration $K_2$. It follows that,
	$$\mathbb{P}[|f(x_0^t)- f(x^*)|>\delta \ \ i.o.  ]=0. $$ 
	Therefore, $f(x_0^t)- f(x^*)\rightarrow  0 $ w.p.1 as $t\rightarrow \infty$.
	
	\subsection{Proof that $\bar{y}^*_t- f(x^*)\rightarrow  0 $ w.p.1 as $t\rightarrow \infty$}
	$$\{|\bar{y}^*_t- f(x^*)|>2\delta\}\subset \{|\bar{y}^*_t-f(x_0^t)|>\delta\} \bigcup \{ |f(x_0^t)- f(x^*)|>\delta\}  $$
	$$ \cup_{t=n}^\infty\{|\bar{y}^*_t- f(x^*)|>2\delta\}\subset\left\{\cup_{t=n}^\infty \{|\bar{y}^*_t-f(x_0^t)|>\delta\} \right\} \bigcup\left\{ \cup_{t=n}^\infty \{ |f(x_0^t)- f(x^*)|>\delta\}\right\}.  $$
	From D.1, we get,
	$$\lim_{n\rightarrow \infty}  \mathbb{P}[\cup_{t=n}^\infty \{|\bar{y}^*_t-f(x_0^t)|>\delta\} ]=0. $$
	From D.2, we get,
	$$\lim_{n\rightarrow \infty}  \mathbb{P}[\cup_{t=n}^\infty \{ |f(x_0^t)- f(x^*)|>\delta\} ]=0. $$
	Therefore,
	$$\lim_{n\rightarrow \infty}\mathbb{P}[\cup_{t=n}^\infty\{|\bar{y}^*_t- f(x^*)|>2\delta\}]<\lim_{n\rightarrow \infty}  \mathbb{P}[\cup_{t=n}^\infty \{|\bar{y}^*_t-f(x_0^t)|>\delta\} ]+ \lim_{n\rightarrow \infty}  \mathbb{P}[\cup_{t=n}^\infty \{ |f(x_0^t)- f(x^*)|>\delta\} ]=0. $$
	It follows that $\bar{y}^*_t- f(x^*)\rightarrow  0 $ w.p.1 as $t\rightarrow \infty$.


	%
	
	
	%
\end{appendix}

\end{document}